\documentclass[preprint,3p]{elsarticle}

\usepackage{amssymb}
\usepackage{amsmath}
\usepackage[all,cmtip]{xy}
\usepackage{amsthm}
\usepackage{color}
\usepackage{enumitem}
\usepackage{tikz}
\usepackage{appendix}
\usepackage[hyphenbreaks]{breakurl}
\usepackage[colorlinks=true]{hyperref}
\input diagxy

\theoremstyle{plain}
  \newtheorem{thm}{Theorem}[section]
  \newtheorem{lem}[thm]{Lemma}
  \newtheorem{prop}[thm]{Proposition}
  \newtheorem{cor}[thm]{Corollary}
\theoremstyle{definition}
  \newtheorem{defn}[thm]{Definition}
  
  \newtheorem{exmp}[thm]{Example}
  \newtheorem{rem}[thm]{Remark}

\DeclareMathAlphabet{\mathcal}{OMS}{cmsy}{m}{n}
\DeclareMathOperator*{\colim}{colim}

\DeclareMathOperator{\Lan}{Lan}
\DeclareMathOperator{\Ran}{Ran}

\DeclareMathOperator{\id}{id}
\DeclareMathOperator{\ob}{ob}

\makeatletter
\def\ps@pprintTitle{%
 \let\@oddhead\@empty
  \let\@evenhead\@empty
  \def\@oddfoot{\vbox{\hsize=\textwidth\footnotesize
  \vskip 8pt
  \copyright 2020. This manuscript version is made available under the CC-BY-NC-ND 4.0 license \url{https://creativecommons.org/licenses/by-nc-nd/4.0/}. The published version is available at \url{https://doi.org/10.1016/j.topol.2020.107524}.\\
  }}%
  \let\@evenfoot\@oddfoot}
\makeatother

\def\oto{{\bfig\morphism<180,0>[\mkern-4mu`\mkern-4mu;]\place(86,0)[\circ]\efig}}
\def\rto{{\bfig\morphism<180,0>[\mkern-4mu`\mkern-4mu;]\place(82,0)[\mapstochar]\efig}}

\newcommand{\da}{\downarrow}
\newcommand{\Da}{\Downarrow}
\newcommand{\ua}{\uparrow}

\newcommand{\lda}{\swarrow}
\newcommand{\rda}{\searrow}

\newcommand{\rat}{\!\rightarrowtail\!}

\newcommand{\bv}{\bigvee}
\newcommand{\bw}{\bigwedge}

\newcommand{\dv}{\dashv}

\newcommand{\nat}{\natural}

\renewcommand{\phi}{\varphi}

\newcommand{\lam}{\lambda}

\newcommand{\si}{\sigma}

\newcommand{\CA}{\mathcal{A}}
\newcommand{\CB}{\mathcal{B}}
\newcommand{\CC}{\mathcal{C}}
\newcommand{\CD}{\mathcal{D}}
\newcommand{\CE}{\mathcal{E}}

\newcommand{\CQ}{\mathcal{Q}}

\newcommand{\CX}{\mathcal{X}}

\newcommand{\BC}{{\bf C}}
\newcommand{\BD}{{\bf D}}
\newcommand{\BE}{{\bf E}}

\newcommand{\Bf}{{\bf f}}
\newcommand{\Bg}{{\bf g}}
\newcommand{\Bh}{{\bf h}}

\newcommand{\sM}{{\sf M}}
\newcommand{\sP}{{\sf P}}

\newcommand{\sQ}{{\sf Q}}
\newcommand{\sQd}{\sQ^{\dag}}

\newcommand{\sY}{{\sf Y}}

\newcommand{\bbR}{\mathbb{R}}

\newcommand{\App}{{\bf App}}
\newcommand{\Cat}{{\bf Cat}}
\newcommand{\CAT}{{\bf CAT}}
\newcommand{\CATB}{{\CAT\Da_c\CB}}

\newcommand{\DIST}{{\bf DIST}}

\newcommand{\Eqv}{{\bf Eqv}}

\newcommand{\Met}{{\bf Met}}
\newcommand{\Ord}{{\bf Ord}}

\newcommand{\REL}{{\bf REL}}
\newcommand{\Rlt}{{\bf Rlt}}
\newcommand{\Set}{{\bf Set}}

\newcommand{\Sup}{{\bf Sup}}

\newcommand{\Top}{{\bf Top}}
\newcommand{\QCAT}{\CQ\text{-}\CAT}
\newcommand{\QBCAT}{\QB\text{-}\CAT}

\newcommand{\sQCat}{\sQ\text{-}\Cat}

\newcommand{\QDIST}{\CQ\text{-}\DIST}
\newcommand{\QREL}{\CQ\text{-}\REL}
\newcommand{\Vect}{{\bf Vect}}

\newcommand{\co}{{\rm co}}

\newcommand{\op}{{\rm op}}

\newcommand{\dPhi}{\Phi^{\da}}
\newcommand{\uPhi}{\Phi_{\ua}}

\newcommand{\PD}{\sP\CD}
\newcommand{\PE}{\sP\CE}

\newcommand{\PdE}{\sP^{\dag}\CE}

\newcommand{\sYd}{\sY^{\dag}}

\newcommand{\QB}{\CQ_{\CB}}
\newcommand{\Fix}{{\sf Fix}}

\newcommand{\oB}{\overline{\CB}}
\newcommand{\oE}{\overline{\CE}}
\newcommand{\oD}{\overline{\CD}}
\newcommand{\oF}{\overline{F}}

\newcommand{\osi}{\overline{\si}}
\newcommand{\otau}{\overline{\tau}}
\newcommand{\omu}{\overline{\mu}}
\newcommand{\olam}{\overline{\lam}}
\newcommand{\ME}{\sM\CE}
\newcommand{\MPhi}{\sM\Phi}

\newcommand{\with}{\mathrel{\&}}

\numberwithin{equation}{section}

\begin{document}

\begin{frontmatter}



\title{Density in categorical topology via quantaloid-enriched categories}


\author{Hongliang Lai}
\ead{hllai@scu.edu.cn}

\author{Lili Shen\corref{cor}}
\ead{shenlili@scu.edu.cn}

\cortext[cor]{Corresponding author.}
\address{School of Mathematics, Sichuan University, Chengdu 610064, China}

\begin{abstract}
Based on Garner's discovery that topological categories are total categories enriched in a quantaloid, this paper presents a series of results related to initial and final density in categorical topology via (co)density in quantaloid-enriched categories, focusing on (co-)Sierpi{\'n}ski objects, Galois correspondences and their fixed points.
\end{abstract}


\end{frontmatter}




\section{Introduction}

Garner's discovery \cite{Garner2014} that topological categories are total categories enriched in a quantaloid has brought new light into the field of \emph{categorical topology} \cite{Adamek1990,Lowen-Colebunders2001} initiated by Br{\"u}mmer \cite{Brummer1971},  Wyler \cite{Wyler1971,Wyler1971a} and Herrlich \cite{Herrlich1974}. As depicted by Garner \cite{Garner2014}, concrete categories over a base category $\CB$ may be interpreted as ``preorders relative to $\CB$'', and topological categories over $\CB$ then become ``complete preorders relative to $\CB$''. To formalize this idea, the theory of \emph{enriched categories} \cite{Eilenberg1966,Kelly1982} and the notion of \emph{totality} \cite{Street1978} are employed in \cite{Garner2014}, and it is shown that concrete categories over $\CB$ may be seen as categories enriched in the \emph{free quantaloid} \cite{Rosenthal1991} $\QB$ induced by $\CB$, while the topologicity of a concrete category over $\CB$ corresponds precisely to the totality of categories enriched in $\QB$ (see \cite{Rosenthal1996,Stubbe2005,Stubbe2006} for the theory of quantaloid-enriched categories). Garner's discovery, as well as being an elegant result itself, unlocks the powerful arsenal of enriched category theory in the study of categorical topology. Following this way, the paper \cite{Shen2016} of Shen and Tholen expands upon Garner's result and interprets several key results of categorical topology from the viewpoint of quantaloid-enriched category theory.

The purpose of this paper is to continue the exploration on categorical topology through the quantaloidal approach, and we focus on the notion of \emph{density} in both contexts. Recall that:
\begin{itemize}
\item a full subcategory $\CD$ of a concrete category $\CE$ over $\CB$ is \emph{initially dense} \cite{Adamek1990} if every $\CE$-object is the domain of some initial source with codomains in $\CD$, and \emph{final density} is the dual concept;
\item a functor $F:\CD\to\CE$ enriched in a quantaloid $\CQ$ is \emph{dense} \cite{Shen2013a,Garner2014,Shen2016} if every $\CE$-object is the colimit of $F$ weighted by a presheaf on $\CD$, and \emph{codensity} is the dual concept.
\end{itemize}
As implicitly suggested in \cite{Garner2014,Shen2016}, initial density and final density in concrete categories amount to codensity and density in quantaloid-enriched categories, respectively. Based on the fundamental 2-equivalence between concrete categories over $\CB$ and categories enriched in the quantaloid $\QB$ elaborated in Section \ref{Concrete-cateogries}, we define \emph{sources} and \emph{sinks} on a general quantaloid-enriched category in Section \ref{Sources-sinks}, which are actually the \emph{discrete} version of (co)presheaves introduced by Stubbe \cite{Stubbe2005}. As revealed in Section \ref{Density}, for a concrete category $\CE$ over $\CB$, sources (resp. sinks) on the corresponding $\QB$-category $\oE$ just disguise $\CE$-structured sources (resp. sinks) in $\CB$. If we slightly generalize the notions of initial and final density to concrete functors (Definition \ref{initial-dense-functor-def}), then it can be shown that a concrete functor $F:\CD\to\CE$ over $\CB$ is initially (resp. finally) dense if, and only if, the corresponding $\QB$-functor $\oF:\oD\to\oE$ is codense (resp. dense); see Proposition \ref{initially-dense-codense}.

Connecting the density in categorical topology and enriched category theory allows us to reconcile several related notions and results in the two lines of research. It is well known that the Sierpi{\'n}ski space, as a one-object full subcategory of the category $\Top$ of topological spaces and continuous maps, is initially dense. In general, a \emph{Sierpi{\'n}ski object} \cite{Manes1976} in a concrete category $\CE$ over $\CB$ is an $\CE$-object that is initially dense in $\CE$ as a one-object full subcategory. With Proposition \ref{initially-dense-codense} at hand, (co-)Sierpi{\'n}ski objects may be defined for categories enriched in any quantaloid (Definition \ref{QCat-Sierpinski-def}). We characterize (co-)Sierpi{\'n}ski objects in quantaloid-enriched categories via (co)tensors (Proposition \ref{Sierpinski-cotensor}) which, conversely, corresponds to an obvious description of (co-)Sierpi{\'n}ski objects in concrete categories (Proposition \ref{Sierpinski-initial-source}).

\emph{Galois correspondences} \cite{Adamek1990} between concrete categories are the concrete version of adjoint functors, which correspond precisely to adjoint functors enriched in a quantaloid. For any quantaloid $\CQ$, it can be shown that left and right adjoint $\CQ$-functors are completely determined by their values on dense and codense $\CQ$-subcategories of their domains, respectively (Proposition \ref{left-adjoint-dense}). In the context of concrete categories, it in particular means that the number of Galois correspondences between certain concrete categories is limited by their actions on a (co-)Sierpi{\'n}ski object, whenever it exists. As applications, we show that there are exactly two Galois correspondences on $\Top$ (Example \ref{Galois-cor-Top}) and three Galois correspondences on the category $\Ord$ of (pre)ordered sets and monotone maps (Example \ref{Galois-cor-Ord}).

Finally, based on our recent paper \cite{Lai2017}, in Section \ref{Fixed-points} we develop representation theorems for fixed points of Galois correspondences between topological categories through initially dense and finally dense concrete functors (Theorems \ref{concrete-representation} and \ref{concrete-representation-dense}), and illustrate them with Examples \ref{Top-Top-cor}, \ref{Ord-Ord-cor} and \ref{Ord-Top-cor}. As an application, the characterization of the MacNeille completion of a concrete category (cf. \cite{Herrlich1976,Adamek1990,Garner2014}) is reproduced as a direct consequence of our general representation theorems (Corollary \ref{ME-representation-concrete}).

For the size issues encountered in this paper, we choose to follow the treatment of \cite{Shen2016}; that is, we employ no specific strict set-theoretical regime, distinguishing only between sets (``small'') and classes (``possibly large''). If a category turns out to be illegitimately large, we add the prefix ``meta'' to indicate that it has to be placed in a higher universe.

\section{Concrete categories as quantaloid-enriched categories} \label{Concrete-cateogries}

A \emph{quantaloid} \cite{Rosenthal1996} is a category enriched in the monoidal closed category $\Sup$ \cite{Joyal1984} of complete lattices and join-preserving functions. Explicitly, a quantaloid $\CQ$ is a (possibly large) 2-category with its 2-cells given by order, such that each hom-set $\CQ(S,T)$ is a complete lattice and the composition of arrows preserves joins in each variable. The induced right adjoints of the composition functions
\begin{align*}
(-\circ u)\dv(-\lda u):&\ \CQ(S,U)\to\CQ(T,U),\\
(v\circ -)\dv(v\rda -):&\ \CQ(S,U)\to\CQ(S,T)
\end{align*}
satisfy
\begin{equation} \label{Q-arrow-adjoint}
v\circ u\leq w\iff v\leq w\lda u\iff u\leq v\rda w
\end{equation}
for all $\CQ$-arrows $u:S\to T$, $v:T\to U$, $w:S\to U$, where $\lda$ and $\rda$ are called \emph{left} and \emph{right implications} of $\CQ$, respectively.

Given a quantaloid $\CQ$, a \emph{$\CQ$-category} (i.e., a \emph{category enriched in $\CQ$}) $\CE$ is given by
\begin{itemize}
\item a class $\CE_0$ of objects,
\item a function $|\text{-}|:\CE_0\to\CQ_0(=\ob\CQ)$ sending each $X\in\CE_0$ to its \emph{extent} $|X|\in\CQ_0$,
\item a family $\CE(X,Y)\in\CQ(|X|,|Y|)$ $(X,Y\in\CE_0)$ of morphisms satisfying
    $$1_{|X|}\leq\CE(X,X)\quad\text{and}\quad\CE(Y,Z)\circ\CE(X,Y)\leq\CE(X,Z)$$
    for all $X,Y,Z\in\CE_0$.
\end{itemize}
It is well known that every (locally small) category $\CB$ induces a \emph{free quantaloid} \cite{Rosenthal1991} $\QB$ given by
\begin{itemize}
\item $\ob\QB=\CB_0(=\ob\CB)$,
\item $\QB(S,T)=\{\Bf\mid\Bf\subseteq\CB(S,T)\}$ for all $S,T\in\CB_0$,
\item $\Bg\circ\Bf=\{g\circ f\mid f\in\Bf,\ g\in\Bg \}$ and ${\bf 1}_S=\{1_S\}$ for all $\Bf\in\QB(S,T)$, $\Bg\in\QB(T,U)$.
\end{itemize}
Considering $\CB$ as a base category, $\QB$-categories are actually \emph{concrete categories over $\CB$}. Explicitly, since a concrete category $\CE$ over $\CB$ \cite{Adamek1990} is given by a faithful functor
$$|\text{-}|:\CE\to\CB,$$
if we identify $\CE(X,Y)$ with $|\CE(X,Y)|\subseteq\CB(|X|,|Y|)$, then $\CE$ is fully described by
\begin{itemize}
\item a class $\CE_0$ of objects,
\item a function $|\text{-}|:\CE_0\to\CB_0$ sending each $X\in\CE_0$ to its \emph{extent} $|X|\in\CB_0$,
\item a family $\CE(X,Y)\subseteq\CB(|X|,|Y|)$ $(X,Y\in\CE_0)$ of morphisms satisfying
    $$1_{|X|}\in\CE(X,X)\quad\text{and}\quad\CE(Y,Z)\circ\CE(X,Y)\subseteq\CE(X,Z)$$
    for all $X,Y,Z\in\CE_0$.
\end{itemize}
We refer to $\CB$-arrows as \emph{maps} and $\CE$-arrows as \emph{morphisms}, respectively; hence, a map $f\in\CB(|X|,|Y|)$ is a morphism (or more precisely, an \emph{$\CE$-morphism}) if $f\in\CE(X,Y)$.

Given a quantaloid $\CQ$, every $\CQ$-category $\CE$ admits an underlying (pre)order on $\CE_0$ given by
\begin{equation} \label{QCat-underlying-order}
X\leq Y\iff |X|=|Y|\quad\text{and}\quad 1_{|X|}\leq\CE(X,Y).
\end{equation}
We write $X\cong Y$ if $X\leq Y$ and $Y\leq X$, and $\CE$ is called \emph{separated} if $X=Y$ whenever $X\cong Y$ in $\CE$.

A \emph{$\CQ$-functor} $F:\CD\to\CE$ between $\CQ$-categories is a function $F:\CD_0\to\CE_0$ such that
$$|X|=|FX|\quad\text{and}\quad\CD(X,Y)\leq\CE(FX,FY)$$
for all $X,Y\in\CD_0$. With the pointwise order of $\CQ$-functors given by
\begin{equation} \label{QCat-2cell}
F\leq G:\CD\to\CE\iff\forall X\in\CD_0:\ FX\leq GX\iff\forall X\in\CD_0:\ 1_{|X|}\leq\CE(FX,GX),
\end{equation}
we obtain a 2-(meta)category
$$\QCAT.$$
With $\CQ=\QB$ for a base category $\CB$, $\QB$-functors precisely characterize \emph{concrete functors} $F:\CD\to\CE$ between concrete categories over $\CB$; that is, functions $F:\CD_0\to\CE_0$ with
$$|X|=|FX|\quad\text{and}\quad\CD(X,Y)\subseteq\CE(FX,FY)$$
for all $X,Y\in\CD_0$. In fact, let
$$\CATB$$
denote the 2-(meta)category of concrete categories and concrete functors over $\CB$ whose 2-cells are given by trading $1_{|X|}\leq\CE(FX,GX)$ with $1_{|X|}\in\CE(FX,GX)$ in \eqref{QCat-2cell}, we obtain a 2-equivalence:

\begin{prop} {\rm\cite{Garner2014}} \label{CATB-equiv-QBCAT}
The 2-(meta)categories $\CATB$ and $\QBCAT$ are 2-equivalent.
\end{prop}

In order to avoid confusion, from now on we write $\oE$ for the corresponding $\QB$-category of each concrete category $\CE$ over $\CB$, and write
$$\oF:\oD\to\oE$$
for the corresponding $\QB$-functor of each concrete functor $F:\CD\to\CE$.

\section{Sources and sinks on quantaloid-enriched categories} \label{Sources-sinks}

From every quantaloid $\CQ$ we may formulate a (meta)quantaloid\footnote{As indicated at the end of the introduction, a \emph{(meta)quantaloid} refers to a quantaloid that lives in a higher universe. For a (meta)quantaloid $\CQ$, each $\CQ(S,T)$ $(S,T\in\CQ_0)$ is an ordered (possibly large) class in which the supremum (or equivalently, infimum) of any subclass exists.}
$$\QREL$$
consisting of the following data:
\begin{itemize}
\item objects of $\QREL$ are classes over $\CQ_0$, i.e., classes $\BE$ equipped with a function $|\text{-}|:\BE\to\CQ_0$;
\item an arrow $\Phi:\BD\rto\BE$ in $\QREL$ is a \emph{$\CQ$-relation}, i.e., a family of $\CQ$-arrows
    $$\Phi(X,Y):|X|\to|Y|\quad(X\in\BD,\ Y\in\BE);$$
\item compositions and implications in $\QREL$ are given by
    \begin{align*}
    &\Psi\circ\Phi:\BC\rto\BE,\quad(\Psi\circ\Phi)(X,Z)=\bv_{Y\in\BD}\Psi(Y,Z)\circ\Phi(X,Y),\\
    &\Xi\lda\Phi:\BD\rto\BE,\quad(\Xi\lda\Phi)(Y,Z)=\bw_{X\in\BC}\Xi(X,Z)\lda\Phi(X,Y),\\
    &\Psi\rda\Xi:\BC\rto\BD,\quad (\Psi\rda\Xi)(X,Y)=\bw_{Z\in\BE}\Psi(Y,Z)\rda\Xi(X,Z)
    \end{align*}
    for all $\CQ$-relations $\Phi:\BC\rto\BD$, $\Psi:\BD\rto\BE$, $\Xi:\BC\rto\BE$;
\item each $\QREL(\BD,\BE)$ is equipped with the order inherited from $\CQ$, i.e.,
    $$\Phi\leq\Psi:\BD\rto\BE\iff\forall X\in\BD,\ \forall Y\in\BE:\ \Phi(X,Y)\leq\Psi(X,Y);$$
\item the identity $\CQ$-relation on $\BE$ is given by
    $$\id_{\BE}:\BE\rto\BE,\quad \id_{\BE}(X,Y)=\begin{cases}
    1_{|X|} & \text{if}\ X=Y,\\
    \bot & \text{else}.
    \end{cases}$$
\end{itemize}
Using the language of $\CQ$-relations, a $\CQ$-category $\CE$ is actually given by a $\CQ$-relation $\CE:\CE_0\rto\CE_0$ with
\begin{equation} \label{Q-cat-def-rel}
\id_{\CE_0}\leq\CE\quad\text{and}\quad\CE\circ\CE\leq\CE.
\end{equation}
Given $\CQ$-categories $\CD$, $\CE$, a $\CQ$-relation $\Phi:\CD_0\rto\CE_0$ becomes a \emph{$\CQ$-distributor} $\Phi:\CD\oto\CE$ if
\begin{equation} \label{Q-dist-def-rel}
\CE\circ\Phi\circ\CD\leq\Phi.
\end{equation}
$\CQ$-categories and $\CQ$-distributors constitute a (meta)quantaloid
$$\QDIST,$$
in which the compositions and implications are calculated in the same way as in $\QREL$, and the identity $\CQ$-distributor on $\CE$ is given by its hom
$$\CE:\CE\oto\CE.$$
Since each class $\BE$ over $\CQ_0$ may be viewed as a \emph{discrete} $\CQ$-category with the hom given by $\id_{\BE}$, $\QREL$ is embedded into $\QDIST$ as a full subquantaloid.

With \eqref{Q-cat-def-rel} and \eqref{Q-dist-def-rel} it is easy to observe the following fact:

\begin{prop} \label{Q-rel-is-dist}
Let $\CD$, $\CE$ be $\CQ$-categories. For any $\CQ$-relation $\Phi:\CD_0\rto\CE_0$,
$$\CE\circ\Phi\circ\CD:\CD\oto\CE$$
is a $\CQ$-distributor.
\end{prop}

For each object $T\in\CQ_0$, we denote by $\{T\}$ the one-object $\CQ$-category whose only object has extent $T$ and hom $1_T$.

\begin{defn} \label{QCat-source-sink-def}
A \emph{sink} (resp. \emph{source}) of extent $T$ on a $\CQ$-category $\CE$ is a $\CQ$-relation
$$\si:\CE_0\rto\{T\}\quad(\text{resp.}\ \ \tau:\{T\}\rto\CE_0).$$
The \emph{supremum} (resp. \emph{infimum}) of $\si$ (resp. $\tau$), when it exists, is an object $\sup\si\in\CE_0$ (resp. $\inf\tau\in\CE_0$) of extent $T$, such that
\begin{equation} \label{sup-inf-def}
\CE(\sup\si,-)=\CE\lda\si\quad (\text{resp.}\ \ \CE(-,\inf\tau)=\tau\rda\CE).
\end{equation}
\end{defn}

A sink $\si$ (resp. source $\tau$) on a $\CQ$-category $\CE$ is called a \emph{presheaf} (resp. \emph{copresheaf}) on $\CE$ if it is a $\CQ$-distributor. It follows from Proposition \ref{Q-rel-is-dist} that every sink $\si:\CE_0\rto\{T\}$ (resp. source $\tau:\{T\}\rto\CE_0$) on $\CE$ generates a presheaf (resp. copresheaf)
$$\si\circ\CE:\CE\oto\{T\}\quad(\text{resp.}\ \ \CE\circ\tau:\{T\}\oto\CE).$$
Suprema (resp. Infima) of sinks (resp. sources) are certainly compatible with suprema (resp. infima) of presheaves (resp. copresheaves) defined in \cite{Shen2013a}:

\begin{prop} \label{sup-si=sup-E-si}
For any sink $\si$ (resp. source $\tau$) on a $\CQ$-category $\CE$,
$$\sup\si=\sup(\si\circ\CE)\quad (\text{resp.}\ \ \inf\tau=\inf(\CE\circ\tau)).$$
\end{prop}

\begin{proof}
Since
$$\CE\lda\si=(\CE\lda\CE)\lda\si=\CE\lda(\si\circ\CE)\quad\text{and}\quad\tau\rda\CE=\tau\rda(\CE\rda\CE)=(\CE\circ\tau)\rda\CE,$$
the conclusion follows soon from \eqref{sup-inf-def}.
\end{proof}

Given a $\CQ$-category $\CE$, it is well known \cite{Stubbe2005} that every presheaf on $\CE$ has a supremum if, and only if, every copresheaf on $\CE$ has an infimum. Following the terminology of \cite{Street1978,Garner2014}, we say that:

\begin{defn} \label{QCat-total-def}
A $\CQ$-category $\CE$ is \emph{total} if every sink on $\CE$ has a supremum, or equivalently, if every source on $\CE$ has an infimum.
\end{defn}

Every $\CQ$-functor $F:\CD\to\CE$ gives rise to a pair of $\CQ$-distributors
$$F_{\nat}:\CD\oto\CE,\quad F_{\nat}(X,Y)=\CE(FX,Y)\quad\text{and}\quad F^{\nat}:\CE\oto\CD,\quad F^{\nat}(Y,X)=\CE(Y,FX),$$
called the \emph{graph} and \emph{cograph} of $F$, respectively. It is easy to see that
\begin{equation} \label{F-leq-G-graph}
F\leq G:\CD\to\CE\iff G_{\nat}\leq F_{\nat}:\CD\oto\CE\iff F^{\nat}\leq G^{\nat}:\CE\oto\CD,
\end{equation}
and thus there are 2-functors
$$(-)_{\nat}:(\QCAT)^{\co}\to\CQ\text{-}\DIST\quad\text{and}\quad(-)^{\nat}:(\QCAT)^{\op}\to\CQ\text{-}\DIST,$$
which are both neutral on objects. Here ``$\co$'' refers to the dualization of 2-cells. Since
$$\CD\leq F^{\nat}\circ F_{\nat}\quad\text{and}\quad F_{\nat}\circ F^{\nat}\leq\CE,$$
we have an adjunction $F_{\nat}\dv F^{\nat}$ in $\CQ\text{-}\DIST$.

\begin{defn} \label{colim-lim-sink-source-def}
The \emph{colimit} (resp. \emph{limit}) of a $\CQ$-functor $F:\CD\to\CE$ weighted by a sink $\si:\CD_0\rto\{T\}$ (resp. source $\tau:\{T\}\rto\CD_0$) on $\CD$, when it exists, is an object $\colim_{\si}F\in\CE_0$ (resp. $\lim_{\tau}F\in\CE_0$) of extent $T$, such that
\begin{equation} \label{colim-lim-def}
\CE({\colim}_{\si}F,-)=F_{\nat}\lda\si\quad(\text{resp.}\ \ \CE(-,{\lim}_{\tau}F)=\tau\rda F^{\nat}).
\end{equation}
\end{defn}

Analogously to Proposition \ref{sup-si=sup-E-si}, we may prove that Definition \ref{colim-lim-sink-source-def} extends the notion of weighted (co)limit in \cite{Stubbe2005}:

\begin{prop} \label{colim-si-F=colim-si-D-F}
For any $\CQ$-functor $F:\CD\to\CE$ and any sink $\si$ (resp. source $\tau$) on $\CD$,
$${\colim}_{\si}F={\colim}_{\si\circ\CD}F\quad (\text{resp.}\ \ {\lim}_{\tau}F={\lim}_{\CD\circ\tau}F).$$
\end{prop}

\begin{proof}
Since
$$F_{\nat}\lda\si=(F_{\nat}\lda\CD)\lda\si=F_{\nat}\lda(\si\circ\CD)\quad\text{and}\quad\tau\rda F^{\nat}=\tau\rda(\CD\rda F^{\nat})=(\CD\circ\tau)\rda F^{\nat},$$
the conclusion follows soon from \eqref{colim-lim-def}.
\end{proof}

It is not difficult to see that weighted colimits (resp. limits) and suprema (resp. infima) can be represented by each other:

\begin{prop} \label{sup-colim-colim-sup}
For each sink $\si$ (resp. source $\tau$) on $\CE$,
\begin{equation} \label{sup-colim}
\sup\si={\colim}_{\si}1_{\CE}\quad(\text{resp.}\ \ \inf\tau={\lim}_{\tau}1_{\CE}).
\end{equation}
Conversely, for each $\CQ$-functor $F:\CD\to\CE$ and sink $\si$ (resp. source $\tau$) on $\CD$,
\begin{equation} \label{colim-sup}
{\colim}_{\si}F=\sup(\si\circ F^{\nat})\quad(\text{resp.}\ \ {\lim}_{\tau}F=\inf(F_{\nat}\circ\tau)).
\end{equation}
\end{prop}

\section{Density in concrete categories and quantaloid-enriched categories} \label{Density}

Given a concrete category $\CE$ over $\CB$, an \emph{$\CE$-structured source} \cite{Adamek1990} $\tau$ is given by a $\CB$-object $T$ and a (possibly large) family
\begin{equation} \label{structured-source}
(g_i:T\to|X_i|)_{i\in I}
\end{equation}
of maps in $\CB$, where $X_i\in\CE_0$ $(i\in I)$. A \emph{lifting} of $\tau$ is an $\CE$-object $Y$ with $|Y|=T$ such that every $g_i$ is an $\CE$-morphism, and the lifting is \emph{initial} if
$$(g_i:Y\to X_i)_{i\in I}$$
is an \emph{initial source} in $\CE$; that is, for any $Z\in\CE_0$, a map $f:|Z|\to|Y|$ becomes an $\CE$-morphism as soon as all maps $g_i\circ f:|Z|\to |X_i|$ are $\CE$-morphisms:
$$\bfig
\ptriangle/->`<-`<-/<800,400>[Y`X_i`Z;g_i`f`g_i\circ f]
\ptriangle(2000,0)/->`<-`<-/<800,400>[|Y|`|X_i|`|Z|;g_i`f`g_i\circ f]
\morphism(1200,200)/|->/<400,0>[`;|\text{-}|]
\efig$$

Note that an $\CE$-structured source $\tau$ given by \eqref{structured-source} may be regrouped into an $\CE_0$-indexed family, which we denote by $\otau$, with
\begin{equation} \label{sturctured-source-obE-index}
\otau(X):=\{g_i:T\to|X_i|\mid i\in I,\ X_i=X\}\subseteq\CB(T,|X|)\quad(X\in\CE_0),
\end{equation}
where $\otau(X)$ is allowed to be empty if no $X_i$ $(i\in I)$ has extent $X$. Thus, $\otau$ is precisely a $\QB$-relation
$$\otau:\{T\}\rto\oE_0;$$
that is, a source on the $\QB$-category $\oE$. From \eqref{Q-arrow-adjoint} it is not difficult to deduce that
\begin{equation} \label{QB-imp-def}
\Bh\lda\Bf=\{g\in\CB(T,U)\mid\forall f\in\Bf:\ g\circ f\in\Bh\}\quad\text{and}\quad\Bg\rda\Bh=\{f\in\CB(S,T)\mid\forall g\in\Bg:\ g\circ f\in\Bh\}.
\end{equation}
for all $\Bf\subseteq\CB(S,T)$, $\Bg\subseteq\CB(T,U)$, $\Bh\subseteq\CB(S,U)$. Hence, the initial lifting of $\tau$, when it exists, is an object $Y\in\CE_0$ satisfying
\begin{align*}
\oE(Z,Y)&=\{f\in\CB(|Z|,|Y|)\mid\forall i\in I:\ g_i\circ f\in\oE(Z,X_i)\}\\
&=\{f\in\CB(|Z|,|Y|)\mid\forall X\in\CE_0,\ \forall g\in\otau(X):\ g\circ f\in\oE(Z,X)\}\\
&=\bigcap\limits_{X\in\CE_0}\{f\in\CB(|Z|,|Y|)\mid\forall g\in\otau(X):\ g\circ f\in\oE(Z,X)\}\\
&=\bigcap\limits_{X\in\CE_0}\otau(X)\rda\oE(Z,X)\\
&=\otau\rda\oE(Z,-)
\end{align*}
for all $Z\in\CE_0$; that is,
\begin{equation} \label{initial-lifting-inf}
\oE(-,Y)=\otau\rda\oE.
\end{equation}
Comparing with Equation \eqref{sup-inf-def} in Definition \ref{QCat-source-sink-def} we have actually proved:

\begin{prop} \label{initial-lifting-inf-prop}
Let $\CE$ be a concrete category over $\CB$. The initial lifting of an $\CE$-structured source $\tau$, when it exists, is precisely the infimum of the source
$$\otau:\{T\}\rto\oE_0$$
on the $\QB$-category $\oE$.
\end{prop}

Dually, every \emph{$\CE$-structured sink} $\si$ given by $(f_i:|X_i|\to T)_{i\in I}$ $(X_i\in\CE_0)$ corresponds to a sink
$$\osi:\oE_0\rto\{T\}$$
with
\begin{equation} \label{sturctured-sink-obE-index}
\osi(X):=\{f_i:|X_i|\to T\mid i\in I,\ X_i=X\}\subseteq\CB(|X|,T)\quad(X\in\CE_0)
\end{equation}
on the $\QB$-category $\oE$. A \emph{lifting} of $\si$ is an $\CE$-object $Y$ with $|Y|=T$ such that every $f_i$ is an $\CE$-morphism, and the lifting is \emph{final} if
$$(f_i:X_i\to Y)_{i\in I}$$
is a \emph{final sink} in $\CE$; that is, for any $Z\in\CE_0$, a map $g:|Y|\to|Z|$ becomes an $\CE$-morphism as soon as all maps $g\circ f_i$ are $\CE$-morphisms. The dual version of Proposition \ref{initial-lifting-inf-prop} reads as follows:

\begin{prop} \label{final-lifting-sup-prop}
Let $\CE$ be a concrete category over $\CB$. The final lifting of an $\CE$-structured sink $\si$, when it exists, is precisely the supremum of the sink
$$\osi:\oE_0\rto\{T\}$$
on the $\QB$-category $\oE$.
\end{prop}

Recall that a concrete category $\CE$ is \emph{topological} \cite{Adamek1990} over $\CB$ if every $\CE$-structured source admits an initial lifting, or equivalently, if every $\CE$-structured sink admits a final lifting. From Definition \ref{QCat-total-def}, Propositions \ref{initial-lifting-inf-prop} and \ref{final-lifting-sup-prop} we have:

\begin{prop} {\rm\cite{Garner2014}} \label{topological=total}
A concrete category $\CE$ over $\CB$ is topological if, and only if, the $\QB$-category $\oE$ is total.
\end{prop}

For a concrete category $\CE$ over $\CB$, a full subcategory $\CD$ of $\CE$ is \emph{initially dense} (resp. \emph{finally dense}) \cite{Adamek1990} if every $Y\in\CE_0$ is the domain (resp. codomain) of some initial source $(g_i:Y\to X_i)_{i\in I}$ (resp. final sink $(f_i:X_i\to Y)_{i\in I})$ in $\CE$ with $X_i\in\CD_0$ for all $i\in I$. This notion may be slightly generated to the following:

\begin{defn} \label{initial-dense-functor-def}
A concrete functor $F:\CD\to\CE$ over $\CB$ is \emph{initially dense} (resp. \emph{finally dense}) if every $Y\in\CE_0$ is the domain (resp. codomain) of some initial source $(g_i:Y\to FX_i)_{i\in I}$ (resp. final sink $(f_i:FX_i\to Y)_{i\in I}$) in $\CE$ with $X_i\in\CD_0$ for all $i\in I$.
\end{defn}

Recall that for any quantaloid $\CQ$, a $\CQ$-functor $F:\CD\to\CE$ is \emph{dense} (resp. \emph{codense}) \cite{Shen2013a,Garner2014,Shen2016,Lai2017} if every object $Y\in\CE_0$ is the colimit (resp. limit) of $F$ weighted by some presheaf (resp. copresheaf) $\tau$ on $\CD$. In particular, a $\CQ$-subcategory $\CD$ of $\CE$ is \emph{dense} (resp. \emph{codense}) if the inclusion $\CQ$-functor $\CD\ \to/^(->/\CE$ is dense (resp. codense).

By Proposition \ref{sup-si=sup-E-si}, it is equivalent to define the (co)density of $\CQ$-functors as follows:

\begin{defn} \label{dense-codense-def}
Let $\CQ$ be a quantaloid. A $\CQ$-functor $F:\CD\to\CE$ is \emph{dense} (resp. \emph{codense}) if every object $Y\in\CE_0$ is the colimit (resp. limit) of $F$ weighted by some sink (resp. source) $\tau$ on $\CD$.
\end{defn}

In order to connect the notions of density in concrete categories and quantaloid-enriched categories, we first prove a lemma:

\begin{lem} \label{source-sink-F}
Let $F:\CD\to\CE$ be a concrete functor over $\CB$. For each $\CE$-structured source $\lam$ (resp. sink $\mu$)
$$(g_i:T\to|FX_i|)_{i\in I}\quad(\text{resp.}\ \ (f_i:|FX_i|\to T)_{i\in I}),$$
where $T\in\CB_0$, $X_i\in\CD_0$ $(i\in I)$, let $\tau$ (resp. $\si$) be the $\CD$-structured source (resp. sink)
$$(g_i:T\to|X_i|)_{i\in I}\quad(\text{resp.}\ \ (f_i:|X_i|\to T)_{i\in I}).$$
Then
$$\oF_{\nat}\circ\otau=\oE\circ\olam\quad(\text{resp.}\ \ \osi\circ\oF^{\nat}=\omu\circ\oE).$$
\end{lem}

\begin{proof}
Note that $|X|=|FX|$ for all $X\in\CD_0$, and thus $\otau(X_i)=\olam(FX_i)$ for all $i\in I$. Moreover, from the definition of $\otau$ and $\olam$ we know that $\otau(X)=\varnothing$ whenever $X\neq X_i$ for any $i\in I$, and $\olam(Y)=\varnothing$ whenever $Y\neq FX_i$ for any $i\in I$. It follows that
\begin{align*}
\oF_{\nat}(-,Y)\circ\otau&=\bigcup_{X\in\CD_0}\oF_{\nat}(X,Y)\circ\otau(X)=\bigcup_{i\in I}\oF_{\nat}(X_i,Y)\circ\otau(X_i)\\
&=\bigcup_{i\in I}\oE(FX_i,Y)\circ\olam(FX_i)=\bigcup_{Z\in\CE_0}\oE(Z,Y)\circ\olam(Z)=\oE(-,Y)\circ\olam
\end{align*}
for all $Y\in\CE_0$, as desired.
\end{proof}

\begin{prop} \label{initially-dense-codense}
A concrete functor $F:\CD\to\CE$ over $\CB$ is initially dense (resp. finally dense) if, and only if, the $\QB$-functor $\oF:\oD\to\oE$ is codense (resp. dense).
\end{prop}

\begin{proof}
If $F:\CD\to\CE$ is initially dense, then for each $Y\in\CE_0$, we have an initial source
$$(g_i:Y\to FX_i)_{i\in I}$$
in $\CE$ with $X_i\in\CD_0$ for all $i\in I$; that is, $Y$ is the initial lifting of the $\CE$-structured source $\lam$ given by
$$(g_i:|Y|\to|FX_i|)_{i\in I}.$$
Let $\tau$ be the $\CD$-structured source
$$(g_i:|Y|\to|X_i|)_{i\in I}.$$
Then
\begin{equation} \label{Y=inf-olam}
Y=\inf\olam=\inf(\oE\circ\olam)=\inf(\oF_{\nat}\circ\otau)={\lim}_{\otau}\oF,
\end{equation}
where the four equalities follow from Propositions \ref{initial-lifting-inf-prop}, \ref{sup-si=sup-E-si}, Lemma \ref{source-sink-F} and Proposition \ref{sup-colim-colim-sup}, respectively. Thus $\oF:\oD\to\oE$ is codense.

Conversely, if $\oF:\oD\to\oE$ is codense, then every $Y\in\CE_0$ satisfies $Y=\lim_{\otau}\oF$ for some source $\otau$ on the $\QB$-category $\oD$, which clearly corresponds to the $\CD$-structured source $\tau$ given by
$$\{g:|Y|\to|X|\mid X\in\CD_0,\ g\in\otau(X)\}.$$
Let $\lam$ be the $\CE$-structured source $\lam$ given by
$$\{g:|Y|\to|FX|\mid X\in\CD_0,\ g\in\otau(X)\}.$$
Then we may deduce analogously to \eqref{Y=inf-olam} that
$$Y={\lim}_{\otau}\oF=\inf(\oF_{\nat}\circ\otau)=\inf(\oE\circ\olam)=\inf\olam.$$
Thus, by Proposition \ref{initial-lifting-inf-prop}, $Y$ is the initial lifting of the $\CE$-structured source $\lam$; that is,
$$\{g:Y\to FX\mid X\in\CD_0,\ g\in\otau(X)\}$$
is an initial source in $\CE$, which completes the proof.
\end{proof}

\section{(Co-)Sierpi{\'n}ski objects of concrete categories via (co)tensors} \label{Sierpinski}

Recall that in a concrete category $\CE$ over $\CB$, a \emph{Sierpi{\'n}ski object} \cite{Manes1976} is an object $X\in\CE_0$ such that $\{X\}$, as a one-object full subcategory of $\CE$, is initially dense in $\CE$. By Proposition \ref{initially-dense-codense}, the notion of Sierpi{\'n}ski object may be defined in quantaloid-enriched categories as follows:

\begin{defn} \label{QCat-Sierpinski-def}
Let $\CQ$ be a quantaloid. In a $\CQ$-category $\CE$, an object $X\in\CE_0$ is \emph{Sierpi{\'n}ski} (resp. \emph{co-Sierpi{\'n}ski}) if $\{X\}$, as a one-object $\CQ$-subcategory of $\CE$, is codense (resp. dense) in $\CE$.
\end{defn}

In order to characterize (co-)Sierpi{\'n}ski objects in $\CQ$-categories, we need the following

\begin{lem} \label{codense-cograph}
A $\CQ$-functor $F:\CD\to\CE$ is dense (resp. codense) if, and only if, $F_{\nat}\lda F_{\nat}=\CD$ (resp. $F^{\nat}\rda F^{\nat}=\CE$).
\end{lem}

\begin{proof}
If $F$ is codense, for each $Y\in\CE_0$ we may find a source $\tau$ on $\CD$ with $Y=\lim_{\tau}F$. Then
\begin{align*}
\CE(-,Y)&\leq F^{\nat}(Y,-)\rda F^{\nat}\\
&\leq\CE(Y,Y)\circ(F^{\nat}(Y,-)\rda F^{\nat})\\
&=(\tau\rda F^{\nat}(Y,-))\circ(F^{\nat}(Y,-)\rda F^{\nat})&(\text{Equation \eqref{colim-lim-def}})\\
&\leq\tau\rda F^{\nat}\\
&=\CE(-,Y),&(\text{Equation \eqref{colim-lim-def}})
\end{align*}
and consequently $\CE(-,Y)=F^{\nat}(Y,-)\rda F^{\nat}=(F^{\nat}\rda F^{\nat})(-,Y)$.

Conversely, if $F^{\nat}\rda F^{\nat}=\CE$, then for each $Y\in\CE_0$ we have
$$\CE(-,Y)=F^{\nat}(Y,-)\rda F^{\nat},$$
which means precisely that $Y=\lim_{F^{\nat}(Y,-)}F$.
\end{proof}

Let $\CE$ be a $\CQ$-category. Recall that the \emph{tensor} (resp. \emph{cotensor}) \cite{Stubbe2006} of a $\CQ$-arrow $f:|X|\to T$ (resp. $g:T\to|X|$) and $X\in\CE_0$, denoted by $f\otimes X$ (resp. $g\rat X$), is an $\CE$-object of extent $T$ such that
\begin{equation} \label{tensor-cotensor-def}
\CE(f\otimes X,-)=\CE(X,-)\lda f\quad(\text{resp.}\ \ \CE(-,g\rat X)=g\rda\CE(-,X)).
\end{equation}

\begin{prop} \label{Sierpinski-cotensor}
Let $\CE$ be a $\CQ$-category. An object $X\in\CE_0$ Sierpi{\'n}ski (resp. co-Sierpi{\'n}ski) if, and only if,
$$Y=\CE(Y,X)\rat X\quad(\text{resp.}\ \ Y=\CE(X,Y)\otimes X)$$
for all $Y\in\CE_0$.
\end{prop}

\begin{proof}
By Lemma \ref{codense-cograph}, $\{X\}$ being codense in $\CE$ precisely translates to
\begin{equation} \label{E=EX-rda-EX}
\CE=\CE(-,X)\rda\CE(-,X);
\end{equation}
that is,
$$\CE(-,Y)=\CE(Y,X)\rda\CE(-,X)$$
for all $Y\in\CE_0$, which means that $Y=\CE(Y,X)\rat X$ by \eqref{tensor-cotensor-def}.
\end{proof}

\begin{exmp} \label{Q-Sierpinski}
Let
$$\sQ=(\sQ,\with,k)$$
be a (unital) \emph{quantale} \cite{Rosenthal1990}, i.e., a one-object quantaloid, where $\with$ is the composition of the unique hom-set on $\sQ$ with $k$ being the unit. Then there are two $\sQ$-categories with $\sQ$ itself as the underlying set, which we denote by $\sQ$ and $\sQd$, respectively, with
$$\sQ(p,q)=q\lda p\quad\text{and}\quad\sQd(p,q)=q\rda p$$
for all $p,q\in\sQ$. We assert that $k$ is a Sierpi{\'n}ski object of $\sQd$ and a co-Sierpi{\'n}ski object of $\sQ$.

In fact, for any $p,q\in\sQ$, the tensor (resp. cotensor) of $p$ and $q$ in $\sQ$ (resp. $\sQd$) is given by
$$p\otimes_{\sQ} q=p\with q\quad(\text{resp.}\ \ p\rat_{\sQd}q=q\with p).$$
It follows that
$$q=(q\lda k)\with k=\sQ(k,q)\otimes_{\sQ}k\quad(\text{resp.}\ \ q=k\with (k\rda q)=\sQd(q,k)\rat_{\sQd}k)$$
for all $q\in\sQ$. Hence, $k$ is a co-Sierpi{\'n}ski object of $\sQ$ (resp. a Sierpi{\'n}ski object of $\sQd$) by Proposition \ref{Sierpinski-cotensor}.
\end{exmp}

\begin{exmp} \label{Q-Sierpinski-co-Sierpinski}
For every quantaloid $\CQ$, there is a $\CQ$-category $\CQ^{\top}$ with $\CQ_0$ as the class of objects and
$$\CQ^{\top}(S,T)=\top_{S,T}:S\to T$$
for all $S,T\in\CQ_0$, where $\top_{S,T}$ refers to the top element of the complete lattice $\CQ(S,T)$. In the $\CQ$-category $\CQ^{\top}$, it follows immediately from Equation \eqref{E=EX-rda-EX} (and its dual) that every object $X\in\CQ_0$ is both Sierpi{\'n}ski and co-Sierpi{\'n}ski.
\end{exmp}

For a concrete category $\CE$ over $\CB$, the \emph{tensor} (resp. \emph{cotensor}) of a $\QB$-arrow $\Bf\subseteq\CB(|X|,T)$ (resp. $\Bg\subseteq\CB(T,|X|)$) and $X\in\CE_0$, when it exists, must satisfy
$$\oE(\Bf\otimes X,Z)=\oE(X,Z)\lda\Bf\quad(\text{resp.}\ \ \oE(Z,\Bg\rat X)=\Bg\rda\oE(Z,X))$$
for all $Z\in\CE_0$; that is, a map $g:|\Bf\otimes X|\to|Z|$ (resp. $f:|Z|\to|\Bg\rat X|$) becomes an $\CE$-morphism as soon as all maps $g\circ f:|X|\to|Z|$ $(f\in\Bf)$ (resp. $g\circ f:|Z|\to|X|$ $(g\in\Bg)$) are $\CE$-morphisms. In other words:

\begin{prop} \label{tensor-final-lifting}
Let $\CE$ be a concrete category over $\CB$. The tensor (resp. cotensor) of a $\QB$-arrow $\Bf\subseteq\CB(|X|,T)$ (resp. $\Bg\subseteq\CB(T,|X|)$) and $X\in\CE_0$, when it exists, is precisely the final (resp. initial) lifting of the structured sink (resp. source) $\Bf\subseteq\CB(|X|,T)$ (resp. $\Bg\subseteq\CB(T,|X|)$).
\end{prop}

Translating Proposition \ref{Sierpinski-cotensor} to concrete categories via Proposition \ref{tensor-final-lifting} produces the following:

\begin{prop} \label{Sierpinski-initial-source}
Let $\CE$ be a concrete category over $\CB$. An object $X\in\CE_0$ is Sierpi{\'n}ski (resp. co-Sierpi{\'n}ski) if, and only if, $\CE(Y,X)$ (resp. $\CE(X,Y)$) is an initial source (resp. a final sink) in $\CE$ for all $Y\in\CE_0$.
\end{prop}

\begin{exmp} \label{QCat-Sierpinski}
If $\sQ$ is a quantale, then the category $\sQCat$ of (small) $\sQ$-categories and $\sQ$-functors is topological over $\Set$ (see \cite[Theorem III.3.1.3]{Hofmann2014}), and $\sQ$ is a Sierpi{\'n}ski object in $\sQCat$. Indeed, for any $\sQ$-category $\CE$, the initiality of the source $\sQCat(\CE,\sQ)$ in $\sQCat$ is easily verified (cf. \cite[Exercise III.1.H]{Hofmann2014}). In particular:
\begin{enumerate}[label=(\arabic*)]
\item \label{QCat-Sierpinski:Ord} If ${\bf 2}=\{0\leq 1\}$ is the two-element Boolean algebra, then $\Ord:={\bf 2}\text{-}\Cat$ is the topological category of (pre)ordered sets and monotone maps over $\Set$. It follows that ${\bf 2}$ is a Sierpi{\'n}ski object of $\Ord$, and moreover, ${\bf 2}$ is also a co-Sierpi{\'n}ski object of $\Ord$. Indeed, for every preordered set $X$, the sink $\Ord({\bf 2},X)$ is final in $\Ord$.
\item \label{QCat-Sierpinski:Met} If $[0,\infty]$ is the Lawvere quantale \cite{Lawvere1973}, then $\Met:=[0,\infty]\text{-}\Cat$ is the topological category of (generalized) metric spaces and non-expansive maps over $\Set$. It follows that $[0,\infty]$ itself, as a metric space, is a Sierpi{\'n}ski object of $\Met$.
\end{enumerate}
\end{exmp}

\begin{exmp} \label{Eqv-Sierpinski}
Let $\Eqv$ denote the full subcategory of $\Ord$ whose objects are sets equipped with an equivalence relation. Then $\Eqv$ has a Sierpi{\'n}ski object given by the set $\{0,1\}$ equipped with the discrete order, i.e., the relation $\{(0,0),(1,1)\}$, and a co-Sierpi{\'n}ski object given by the set $\{0,1\}$ equipped with the indiscrete order, i.e., the relation $\{(0,0),(1,1),(0,1),(1,0)\}$.
\end{exmp}

\begin{exmp} \label{B-Sierpinski}
If $\CB$ is considered as a concrete category over itself, then every object $X\in\CB_0$ is both Sierpi{\'n}ski and co-Sierpi{\'n}ski. In fact, the corresponding $\QB$-category $\oB$ is precisely the $\QB$-category $\QB^{\top}$ given in Example \ref{Q-Sierpinski-co-Sierpinski}.
\end{exmp}

\begin{exmp} \label{Top-Sierpinski}
The \emph{Sierpi{\'n}ski space} $S$, whose underlying set is $\{0,1\}$ and whose open sets are $\{\varnothing,\{1\},\{0,1\}\}$, is clearly a Sierpi{\'n}ski object of the category $\Top$ of topological spaces and continuous maps.
\end{exmp}

\begin{exmp} \label{App-Sierpinski}
\emph{Approach spaces} \cite{Lowen1989b,Lowen2015} may be considered as topological spaces valued in the Lawvere quantale $[0,\infty]$ (cf. \cite[Section III.2.4]{Hofmann2014}). In the topological category $\App$ of approach spaces and contractions over $\Set$, $[0,\infty]$ is itself an approach space and serves as a Sierpi{\'n}ski object of $\App$ (see \cite[Theorem 1.3.19]{Lowen2015}).
\end{exmp}

\begin{exmp} \label{Vect-co-Sierpinski}
The category $\Vect$ of real vector spaces and linear transformations is concrete (but not topological) over $\Set$. The real plane $\bbR^2$ is a co-Sierpi{\'n}ski object of $\Vect$ since, for any real vector space $X$, the sink $\Vect(\bbR^2,X)$ is final in $\Vect$.
\end{exmp}

\begin{exmp} \label{Rlt-Sierpinski} \cite{Herrlich1992}
Let $\Rlt$ denote the topological category over $\Set$ whose objects are sets equipped with an arbitrary relation and whose morphisms are relation-preserving maps. Then $\Rlt$ has a Sierpi{\'n}ski object given by the set $\{0,1\}$ equipped with the relation $\{(0,0),(0,1),(1,0)\}$, and a co-Sierpi{\'n}ski object given by the set $\{0,1\}$ equipped with the relation $\{(0,1)\}$.
\end{exmp}

\section{Galois correspondences} \label{Galois-correspondences}

A \emph{Galois correspondence} \cite{Adamek1990} between concrete categories is precisely an internal adjunction in $\CATB$. In other words, concrete functors $L:\CD\to\CE$ and $R:\CE\to\CD$ form a Galois correspondence $(L,R)$ between $\CD$ and $\CE$, if
\begin{equation} \label{Galois-cor-def}
1_{\CD}\leq RL\quad\text{and}\quad LR\leq 1_{\CE};
\end{equation}
or equivalently, if
\begin{equation} \label{Fnat=Gnat}
\CE(L-,-)=\CD(-,R-).
\end{equation}

Galois correspondences between concrete categories clearly amount to \emph{adjoint functors} between quantaloid-enriched categories. In general, for any quantaloid $\CQ$, $\CQ$-functors $L:\CD\to\CE$ and $R:\CE\to\CD$ form an adjunction in $\QCAT$, denoted by $L\dv R$, if \eqref{Galois-cor-def} or \eqref{Fnat=Gnat} holds. In this case, $L$ is called the \emph{left adjoint} of $R$, and $R$ is called the \emph{right adjoint} of $L$.

The following proposition claims that a left (resp. right) adjoint $\CQ$-functor is completely determined by its action on a dense (resp. codense) $\CQ$-subcategory of its domain:

\begin{prop} \label{left-adjoint-dense}
Let $F:\CC\to\CD$ be a dense (resp. codense) $\CQ$-functor. If $L,L':\CD\to\CE$ (resp. $R,R':\CD\to\CE$) are left (resp. right) adjoint $\CQ$-functors satisfying
$$LF\cong L'F\quad(\text{resp.}\ \ RF\cong R'F),$$
then $L\cong L'$ (resp. $R\cong R'$). In particular, if $\CD$ has a co-Sierpi{\'n}ski (resp. Sierpi{\'n}ski) object $X$ and $LX\cong L'X$ (resp. $RX\cong R'X$), then $L\cong L'$ (resp. $R\cong R'$).
\end{prop}

\begin{proof}
For every $X\in\CD_0$, by the density of $F$ we may find a sink $\si$ on $\CC$ with $X=\colim_{\si}F$. Suppose that $L\dv R:\CE\to\CD$ in $\QCAT$, then
$$\CE(LX,-)=\CD(X,R-)=F_{\nat}(-,R-)\lda\si=\CD(F-,R-)\lda\si=\CE(LF-,-)\lda\si=(LF)_{\nat}\lda\si;$$
that is, $LX=\colim_{\si}LF$. Similarly we may deduce that $L'X=\colim_{\si}L'F$. Thus $L\cong L'$ follows from $LF\cong L'F$.
\end{proof}

\begin{cor} \label{Galois-cor-dense}
Let $F:\CC\to\CD$ be a finally dense (resp. an initially dense) concrete functor over $\CB$. If $(L,R)$ and $(L',R')$ are both Galois correspondences between concrete categories $\CD$ and $\CE$ over $\CB$ satisfying
$$LF\cong L'F\quad(\text{resp.}\ \ RF\cong R'F),$$
then $L\cong L'$ (resp. $R\cong R'$). In particular, if $\CD$ has a co-Sierpi{\'n}ski (resp. Sierpi{\'n}ski) object $X$ and $LX\cong L'X$ (resp. $RX\cong R'X$), then $L\cong L'$ (resp. $R\cong R'$).
\end{cor}

\begin{rem} \label{iso-concrete}
In Corollary \ref{Galois-cor-dense} the symbol ``$\cong$'' is used in the context of concrete categories over $\CB$. Throughout this paper, we make the convention that ``$X\cong Y$'' in a concrete category $\CE$ over $\CB$ always refers to $X\leq Y$ and $Y\leq X$ in the underlying order of the $\QB$-category $\oE$ (see \eqref{QCat-underlying-order}), which is stronger than the isomorphism of objects in an ordinary category.
\end{rem}

As applications of Corollary \ref{Galois-cor-dense}, we claim that there are exactly two Galois correspondences (up to homeomorphism) on $\Top$ and three Galois correspondences (up to isomorphism) on $\Ord$:

\begin{exmp} \label{Galois-cor-Top}
If $(L,R)$ is a Galois correspondence on $\Top$, then $R$ is completely determined by its action on the Sierpi{\'n}ski space $S$; that is, $RS$ can only be the two-element set $\{0,1\}$ equipped with the Sierpi{\'n}ski topology itself, the indiscrete topology and the discrete topology:
\begin{enumerate}[label=(\arabic*)]
\item \label{Galois-cor-Top:Sierpinski} If $RS$ is the Sierpi{\'n}ski space itself, then the generated Galois correspondence $(1_{\Top},1_{\Top})$ consists of the identity functors on $\Top$.
\item \label{Galois-cor-Top:indiscrete} If $RS$ endows $\{0,1\}$ with the indiscrete topology, then the generated Galois correspondence $(L,R)$ is given by $L$ sending each topological space $X$ to the discrete space on $|X|$ and $R$ sending each $X$ to the indiscrete space on $|X|$. Indeed, $(L,R)$ is a Galois correspondence since
    $$\Top(LX,Y)=\Set(|X|,|Y|)=\Top(X,RY)$$
    for all topological spaces $X$, $Y$.
\item \label{Galois-cor-Top:discrete} There is no Galois correspondence $(L,R)$ on $\Top$ with $RS$ being the discrete space. In fact, if such $(L,R)$ exists, then
    \begin{equation} \label{Top-LX-S-X-RS}
    \Top(LX,S)=\Top(X,RS)
    \end{equation}
    for any topological space $X$. Note that continuous maps $LX\to S$ are precisely the characteristic functions of open subsets of $LX$, while continuous maps $X\to RS$ are precisely the characteristic functions of clopen subsets of $X$. Hence, Equation \eqref{Top-LX-S-X-RS} means that the topology of $LX$ should exactly consists of all clopen subsets of $X$. However, the collection of clopen subsets of a topological space $X$ may not form a topology on $|X|$; for instance, when $X=\bbR_l$ is the Sorgenfrey line. This negates the existence of such a Galois correspondence.
\end{enumerate}
\end{exmp}

\begin{exmp} \label{Galois-cor-Ord}
Since the two-element Boolean algebra ${\bf 2}$ is a co-Sierpi{\'n}ski object of $\Ord$, a Galois correspondence $(L,R)$ on $\Ord$ is completely determined by $L{\bf 2}$, which could endow the two-element set $\{0,1\}$ with the canonical order $0\leq 1$, the discrete order and the indiscrete order. The induced three Galois correspondences on $\Ord$ are as follows:
\begin{enumerate}[label=(\arabic*)]
\item \label{Galois-cor-Ord:2} If $L{\bf 2}={\bf 2}$, then the generated Galois correspondence $(1_{\Ord},1_{\Ord})$ consists of the identity functors on $\Ord$.
\item \label{Galois-cor-Ord:discrete} If $L{\bf 2}$ endows $\{0,1\}$ with the discrete order, then the generated Galois correspondence $(L,R)$ is given by $L$ sending each ordered set $X$ to the discrete order ``$=$'' on $|X|$ and $R$ sending each $X$ to the indiscrete order on $|X|$ (i.e., the relation $|X|\times|X|$). Indeed, $(L,R)$ is a Galois correspondence since
    $$\Ord(LX,Y)=\Set(|X|,|Y|)=\Ord(X,RY)$$
    for all ordered sets $X$, $Y$.
\item \label{Galois-cor-Ord:indiscrete} If $L{\bf 2}$ endows $\{0,1\}$ with the indiscrete order, then the generated Galois correspondence $(L,R)$ may be described as follows:
    \begin{itemize}
    \item $LX$ sends each ordered set $X=(X,\leq)$ to the equivalence relation $\leq_L$ on $|X|$ generated by $\leq$, i.e., the intersection of the equivalence relations on $|X|$ that contain $\leq$.
    \item $RX$ sends each ordered set $X=(X,\leq)$ to the equivalence relation $\leq_R=\leq\cap\leq^{\op}$ on $|X|$.
    \end{itemize}
    To see that $(L,R)$ is a Galois correspondence on $\Ord$, for any ordered sets $X$, $Y$ and any map $f:|X|\to|Y|$, we need to show that $f\in\Ord(LX,Y)$ if and only if $f\in\Ord(X,RY)$.

    On one hand, if $f:LX\to Y$ is monotone, then for any $x,y\in X$ with $x\leq y$, we have $x\leq_L y$ and $y\leq_L x$, which implies that $fx\leq fy$ and $fy\leq fx$ in $Y$, i.e., $fx\leq_R fy$; that is, $f:X\to RY$ is monotone.

    On the other hand, suppose that $f:X\to RY$ is monotone. For any $x,y\in X$ with $x\leq_L y$, there exist $z_1,z_2,\dots,z_n\in X$ such that $x=z_1$, $z_n=y$ and $z_i\leq z_{i+1}$ or $z_{i+1}\leq z_i$ ($i=1,\dots,n-1$). It follows that $fz_i\leq_R fz_{i+1}$ for all $i=1,\dots,n-1$, and consequently $fx\leq_R fy$, which necessarily forces $fx\leq fy$ in $Y$. This proves the monotonicity of $f:LX\to Y$.
\end{enumerate}
\end{exmp}

\section{Fixed points of Galois correspondences} \label{Fixed-points}

If $F:\CE\to\CE$ is a $\CQ$-functor on a $\CQ$-category $\CE$ (resp. a concrete functor on a concrete category $\CE$ over $\CB$), there is a $\CQ$-subcategory (resp. full subcategory) of $\CE$, denoted by $\Fix(F)$, whose objects are \emph{fixed points} of $F$, i.e.,
$$\Fix(F)_0:=\{X\in\CE_0\mid FX\cong X\}\footnote{See Remark \ref{iso-concrete} for the meaning of the symbol ``$\cong$'' in the context of concrete categories.}.$$

A $\CQ$-functor (resp. concrete functor) $F:\CD\to\CE$ is called an \emph{equivalence} (resp. a \emph{concrete equivalence}) if there is a $\CQ$-functor (resp. concrete functor) $G:\CE\to\CD$ with
\begin{equation} \label{QCat-equiv-def}
1_{\CD}\cong GF\quad\text{and}\quad FG\cong 1_{\CE}.
\end{equation}
In this case, we write $\CD\simeq\CE$ to denote that $\CD$ and $\CE$ are \emph{equivalent} $\CQ$-categories (resp. \emph{concretely equivalent} categories). Moreover, if the symbols ``$\cong$'' in \eqref{QCat-equiv-def} are replaced by ``$=$'', then $F$ is called an \emph{isomorphism} (resp. a \emph{concrete isomorphism}), and we write $\CD\cong\CE$ to denote that $\CD$ and $\CE$ are \emph{isomorphic} $\CQ$-categories (resp. \emph{concretely isomorphic} categories).

\begin{rem}
It should be pointed out that our definition of \emph{concrete equivalence} here is different from \cite[Remark 5.13]{Adamek1990}. More specifically, a concrete equivalence defined by \eqref{QCat-equiv-def} is stronger than a concrete functor that is an equivalence of (ordinary) categories.
\end{rem}

Assuming the axiom of choice, it is not difficult to see that

\begin{prop} {\rm\cite{Stubbe2005}} \label{con-equiv-def-ff}
A $\CQ$-functor (resp. concrete functor) $F:\CD\to\CE$ is an equivalence (resp. a concrete equivalence) if, and only if,
\begin{enumerate}[label={\rm(\arabic*)}]
\item $F$ is \emph{fully faithful} in the sense that $\CD(X,Y)=\CE(FX,FY)$ for all $X,Y\in\CD_0$,
\item $F$ is \emph{essentially surjective} in the sense that there exists $X\in\CD_0$ with $Y\cong FX$ for all $Y\in\CE_0$.
\end{enumerate}
\end{prop}

Given an adjunction $L\dv R:\CE\to\CD$ in $\QCAT$ (resp. a Galois correspondence $(L,R)$ between concrete categories $\CD$ and $\CE$), it is easy to see that the restrictions of $L$ and $R$,
\begin{align*}
L|_{\Fix(RL)}:\Fix(RL)\to\Fix(LR)\quad\text{and}\quad R|_{\Fix(LR)}:\Fix(LR)\to\Fix(RL),
\end{align*}
establish an equivalence of $\CQ$-categories (resp. a concrete equivalence of concrete categories). Thus, objects in both $\Fix(RL)$ and $\Fix(LR)$ will be referred to as fixed points of the adjoint $\CQ$-functors $L\dv R$ (resp. the Galois correspondence $(L,R)$).

\begin{thm} {\rm\cite{Lai2017}} \label{QCat-representation}
Let $L\dv R:\CE\to\CD$ be an adjunction in $\QCAT$. A $\CQ$-category $\CX$ is equivalent to $\Fix(RL)$ if, and only if, there exist essentially surjective $\CQ$-functors $F:\CD\to\CX$ and $G:\CE\to\CX$ with
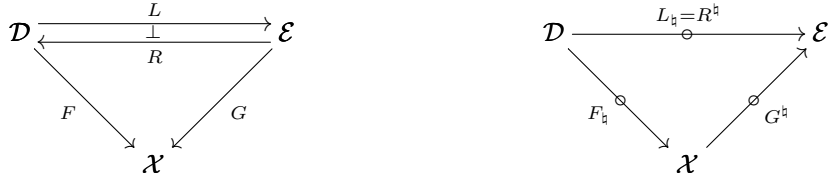
\begin{equation} \label{L=G-F}
L_{\nat}=G^{\nat}\circ F_{\nat},\quad\text{i.e.,}\quad\CE(L-,-)=\CX(F-,G-).
\end{equation}
$$\bfig
\Vtriangle/@{>}@<4pt>`->`->/[\CD`\CE`\CX;L`F`G]
\morphism(1000,500)|b|/@{>}@<3pt>/<-1000,0>[\CE`\CD;R]
\place(500,510)[\mbox{\scriptsize$\bot$}]
\Vtriangle(2000,0)/->`->`<-/[\CD`\CE`\CX;L_{\nat}=R^{\nat}`F_{\nat}`G^{\nat}]
\place(2250,250)[\circ] \place(2750,250)[\circ] \place(2500,500)[\circ]
\efig$$
\end{thm}

\begin{proof}
For the ``only if'' part, let
$$F:\CD\to\Fix(RL)\quad\text{and}\quad G:\CE\to\Fix(RL)$$
be the codomain restriction of $RL:\CD\to\CD$ and $R:\CE\to\CD$, respectively, then $F$ and $G$ are clearly essentially surjective and satisfy the Equation \eqref{L=G-F}. Conversely, for the ``if'' part, the restriction
$$F|_{\Fix(RL)}:\Fix(RL)\to\CX$$
of $F:\CD\to\CX$ is an equivalence of $\CQ$-categories. For details we refer to the proof of \cite[Theorem 3.3]{Lai2017}.
\end{proof}

If $L\dv R:\CE\to\CD$ is an adjunction between total $\CQ$-categories, then $RL:\CE\to\CE$ is a \emph{$\CQ$-closure operator} \cite{Shen2013a} in the sense that
$$1_{\CE}\leq RL\quad\text{and}\quad RLRL\cong RL,$$
and thus $\Fix(RL)$ is also a total $\CQ$-category by \cite[Proposition 3.5]{Shen2013a}. In this case:

\begin{thm} {\rm\cite{Lai2017}} \label{QCat-representation-dense}
Let $L\dv R:\CE\to\CD$ be an adjunction between total $\CQ$-categories. A total $\CQ$-category $\CX$ is equivalent to $\Fix(RL)$ if, and only if, there exist dense $\CQ$-functors $F:\CA\to\CX$, $K:\CA\to\CD$ and codense $\CQ$-functors $G:\CC\to\CX$, $H:\CC\to\CE$ with
\begin{equation} \label{H-L-K=G-F}
H^{\nat}\circ L_{\nat}\circ K_{\nat}=G^{\nat}\circ F_{\nat},\quad\text{i.e.,}\quad\CE(LK-,H-)=\CX(F-,G-).
\end{equation}
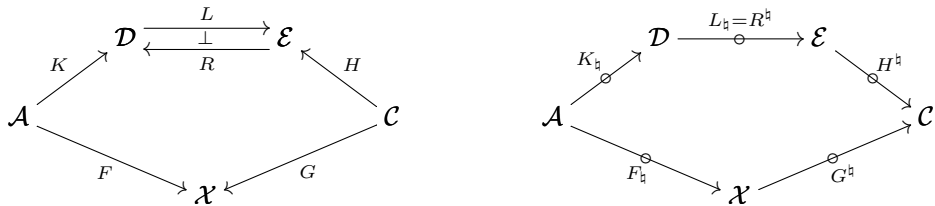
$$\bfig
\morphism(-700,0)<400,300>[\CA`\CD;K]
\morphism(700,0)<-400,300>[\CC`\CE;H]
\morphism(-700,0)|b|<700,-300>[\CA`\CX;F]
\morphism(700,0)|b|<-700,-300>[\CC`\CX;G]
\morphism(-300,300)|a|/@{>}@<4pt>/<600,0>[\CD`\CE;L]
\morphism(300,300)|b|/@{>}@<4pt>/<-600,0>[\CE`\CD;R]
\place(0,305)[\mbox{\scriptsize$\bot$}]
\morphism(1300,0)<400,300>[\CA`\CD;K_{\nat}]
\morphism(1700,300)|a|<600,0>[\CD`\CE;L_{\nat}=R^{\nat}]
\morphism(2300,300)<400,-300>[\CE`\CC;H^{\nat}]
\morphism(1300,0)|b|<700,-300>[\CA`\CX;F_{\nat}]
\morphism(2000,-300)|b|<700,300>[\CX`\CC;G^{\nat}]
\place(1500,150)[\circ] \place(2000,300)[\circ] \place(2500,150)[\circ] \place(1650,-150)[\circ] \place(2350,-150)[\circ]
\efig$$
\end{thm}

\begin{proof}
For the ``only if'' part, we may find essentially surjective $\CQ$-functors $F:\CD\to\CX$ and $G:\CE\to\CX$ with $L_{\nat}=G^{\nat}\circ F_{\nat}$ by Theorem \ref{QCat-representation}. Then the dense $\CQ$-functors $F:\CD\to\CX$, $1_{\CD}:\CD\to\CD$ and the codense $\CQ$-functors $G:\CE\to\CX$, $1_{\CE}:\CE\to\CE$ clearly satisfy \eqref{H-L-K=G-F}. Conversely, for the ``if'' part, let
$$\Lan_K F:\CD\to\CX\quad\text{and}\quad\Ran_H G:\CE\to\CX
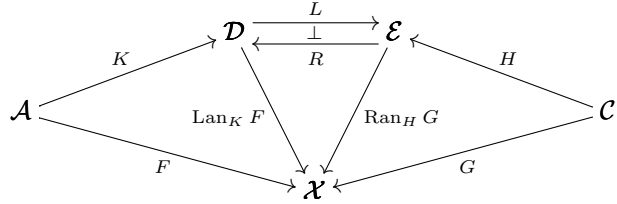
$$
$$\bfig
\morphism(-1100,0)<800,300>[\CA`\CD;K]
\morphism(1100,0)<-800,300>[\CC`\CE;H]
\morphism(-1100,0)|b|<1100,-300>[\CA`\CX;F]
\morphism(1100,0)|b|<-1100,-300>[\CC`\CX;G]
\morphism(-300,300)|a|/@{>}@<4pt>/<600,0>[\CD`\CE;L]
\morphism(300,300)|b|/@{>}@<4pt>/<-600,0>[\CE`\CD;R]
\morphism(-300,300)|l|<300,-600>[\CD`\CX;\Lan_K F]
\morphism(300,300)|r|<-300,-600>[\CE`\CX;\Ran_H G]
\place(0,305)[\mbox{\scriptsize$\bot$}]
\efig$$
denote the \emph{left Kan extension} of $F$ along $K$ and the \emph{right Kan extension} of $G$ along $H$ given by
$$(\Lan_K F)X={\colim}_{K_{\nat}(-,X)}F\quad\text{and}\quad(\Ran_H G)Y={\lim}_{H^{\nat}(Y,-)}G,$$
respectively, whose existence is guaranteed by the totality of $\CX$. Then $\Lan_K F$ and $\Ran_H G$ are both essentially surjective, and satisfy
$$L_{\nat}=(\Ran_H G)^{\nat}\circ(\Lan_K F)_{\nat}.$$
Hence, by Theorem \ref{QCat-representation}, $\CX$ is equivalent to $\Fix(RL)$. For details we refer to the proof of \cite[Theorem 5.1]{Lai2017}.
\end{proof}

By Propositions \ref{CATB-equiv-QBCAT} and \ref{initially-dense-codense}, Theorems \ref{QCat-representation} and \ref{QCat-representation-dense} may be formulated in the context of concrete categories as follows:

\begin{thm} \label{concrete-representation}
Let $(L,R)$ be a Galois correspondence between concrete categories $\CD$ and $\CE$ over $\CB$. A concrete category $\CX$ over $\CB$ is concretely equivalent to $\Fix(RL)$ if, and only if, there exist essentially surjective concrete functors $F:\CD\to\CX$ and $G:\CE\to\CX$ such that
\begin{equation} \label{L=G-F:concrete}
\CE(L-,-)=\CX(F-,G-).
\end{equation}
\end{thm}

\begin{thm} \label{concrete-representation-dense}
Let $(L,R)$ be a Galois correspondence between topological categories $\CD$ and $\CE$ over $\CB$. A topological category $\CX$ over $\CB$ is concretely equivalent to $\Fix(RL)$ if, and only if, there exist finally dense concrete functors $F:\CA\to\CX$, $K:\CA\to\CD$ and initially dense concrete functors $G:\CC\to\CX$, $H:\CC\to\CE$ with
\begin{equation} \label{H-L-K=G-F:concrete}
\CE(LK-,H-)=\CX(F-,G-).
\end{equation}
\end{thm}

\begin{exmp} \label{Top-Top-cor}
For the Galois correspondence $(L,R)$ on $\Top$ given in Example \ref{Galois-cor-Top}\ref{Galois-cor-Top:indiscrete}, $\Fix(RL)$ is concretely equivalent to $\Set$. In fact, since the forgetful functor $|\text{-}|:\Top\to\Set$ is clearly surjective on objects and
$$\Top(LX,Y)=\Set(|X|,|Y|)$$
for all topological spaces $X$, $Y$, the conclusion follows immediately from Theorem \ref{concrete-representation}.
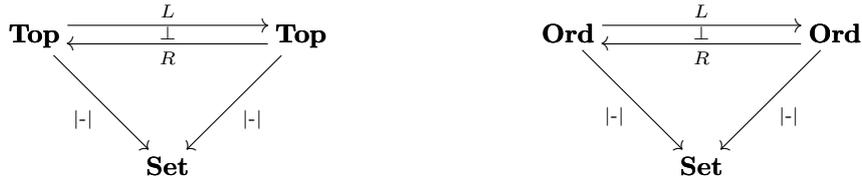
$$\bfig
\Vtriangle/@{>}@<4pt>`->`->/[\Top`\Top`\Set;L`|\text{-}|`|\text{-}|]
\morphism(1000,500)|b|/@{>}@<3pt>/<-1000,0>[\Top`\Top;R]
\place(500,510)[\mbox{\scriptsize$\bot$}]
\Vtriangle(2000,0)/@{>}@<4pt>`->`->/[\Ord`\Ord`\Set;L`|\text{-}|`|\text{-}|]
\morphism(3000,500)|b|/@{>}@<3pt>/<-1000,0>[\Ord`\Ord;R]
\place(2500,510)[\mbox{\scriptsize$\bot$}]
\efig$$
Similarly, for the Galois correspondence $(L,R)$ on $\Ord$ given in Example \ref{Galois-cor-Ord}\ref{Galois-cor-Ord:discrete}, we may also deduce analogously as above that $\Fix(RL)\simeq\Set$.
\end{exmp}

\begin{exmp} \label{Ord-Ord-cor}
For the Galois correspondence $(L,R)$ on $\Ord$ given in Example \ref{Galois-cor-Ord}\ref{Galois-cor-Ord:indiscrete}, $\Fix(RL)$ is concretely equivalent to the category $\Eqv$ given in Example \ref{Eqv-Sierpinski}. In fact, since the two-element Boolean algebra ${\bf 2}$ is simultaneously initially dense and finally dense in $\Ord$ (see Example \ref{QCat-Sierpinski}\ref{QCat-Sierpinski:Ord}), while $L{\bf 2}$ and $R{\bf 2}$ are respectively finally dense and initially dense in $\Eqv$ (see Example \ref{Eqv-Sierpinski}), the conclusion follows from Theorem \ref{concrete-representation-dense} because \eqref{H-L-K=G-F:concrete} reduces to a trivial equation
$$\Ord(L{\bf 2},{\bf 2})=\Eqv(L{\bf 2},R{\bf 2}).$$
$$\bfig
\morphism(-800,0)/^(->/<400,300>[\{\bf 2\}`\Ord;]
\morphism(800,0)/_(->/<-400,300>[\{\bf 2\}`\Ord;]
\morphism(-800,0)|b|<800,-300>[\{\bf 2\}`\Eqv;L|_{\{\bf 2\}}]
\morphism(800,0)|b|<-800,-300>[\{\bf 2\}`\Eqv;R|_{\{\bf 2\}}]
\morphism(-400,300)|a|/@{>}@<4pt>/<800,0>[\Ord`\Ord;L]
\morphism(400,300)|b|/@{>}@<4pt>/<-800,0>[\Ord`\Ord;R]
\place(0,305)[\mbox{\scriptsize$\bot$}]
\efig$$
\end{exmp}

\begin{exmp} \label{Ord-Top-cor}
There is a Galois correspondence $(L,R)$ between $\Ord$ and $\Top$, where $L$ sends each ordered set $X$ to the topological space with all lower subsets of $X$ as its closed subsets, and $R$ sends each topological space $X$ to its specialization order on $|X|$. We claim that $\Fix(RL)$ is concretely equivalent to $\Ord$. In fact, since ${\bf 2}$ is both initially dense and finally dense in $\Ord$, and the Sierpi{\'n}ski space $S$ is initially dense in $\Top$, the conclusion follows from Theorem \ref{concrete-representation-dense} because \eqref{H-L-K=G-F:concrete} now becomes a trivial equation
$$\Top(L{\bf 2},S)=\Ord({\bf 2},RS).$$
$$\bfig
\morphism(-800,0)/^(->/<400,300>[\{\bf 2\}`\Ord;]
\morphism(800,0)/_(->/<-400,300>[\{S\}`\Top;]
\morphism(-800,0)|b|/_(->/<800,-300>[\{\bf 2\}`\Ord;]
\morphism(800,0)|b|<-800,-300>[\{S\}`\Ord;R|_{\{S\}}]
\morphism(-400,300)|a|/@{>}@<4pt>/<800,0>[\Ord`\Top;L]
\morphism(400,300)|b|/@{>}@<4pt>/<-800,0>[\Top`\Ord;R]
\place(0,305)[\mbox{\scriptsize$\bot$}]
\efig$$
\end{exmp}

A topological category $\CD$ over $\CB$ is called the \emph{MacNeille completion} \cite{Herrlich1976,Adamek1990} of its full subcategory $\CE$ if $\CE$ is simultaneously initially dense and finally dense in $\CD$. The MacNeille completion of $\CE$ is unique up to concrete isomorphism (cf. \cite[Exercise 21H]{Adamek1990}), and we denote by
$$\ME$$
the MacNeille completion of $\CE$. In fact, the MacNeille completion is a special case of a much more general construction in quantaloid-enriched categories, as we will describe below.

Given a quantaloid $\CQ$ and a $\CQ$-category $\CE$, presheaves (resp. copresheaves) on $\CE$ constitute a $\CQ$-(meta)category $\PE$ (resp. $\PdE$), with
$$\PE(\si,\si')=\si'\lda\si\quad(\text{resp.}\ \ \PdE(\tau,\tau')=\tau'\rda\tau)$$
for all presheaves $\si,\si'$ (resp. copresheaves $\tau,\tau'$) on $\CE$. The \emph{Yoneda embedding} (resp. \emph{co-Yoneda embedding})
$$\sY_{\CE}:\CE\to\PE,\quad X\mapsto\CE(-,X)\quad(\text{resp.}\ \ \sYd_{\CE}:\CE\to\PdE,\quad X\mapsto\CE(X,-))$$
is a dense (resp. codense) $\CQ$-functor since
$$\si={\colim}_{\si}\sY_{\CE}\quad(\text{resp.}\ \ \tau={\lim}_{\tau}\sYd_{\CE})$$
for all presheaves $\si$ (resp. copresheaves $\tau$) on $\CE$. It is well known that both $\PE$ and $\PdE$ are total (see \cite[Proposition 6.4]{Stubbe2005}).

Each $\CQ$-distributor $\Phi:\CD\oto\CE$ induces an \emph{Isbell adjunction} \cite{Shen2013a}
$$\uPhi\dv\dPhi:\PdE\to\PD$$
given by
$$\uPhi\si=\Phi\lda\si\quad\text{and}\quad\dPhi\tau=\tau\rda\Phi$$
for all presheaves $\si$ on $\CD$ and copresheaves $\tau$ on $\CE$, whose fixed points constitute a total $\CQ$-(meta)category
$$\MPhi:=\Fix(\dPhi\uPhi);$$
explicitly, for any presheaf $\si$ on $\CD$,
$$\si\in(\MPhi)_0\iff\dPhi\uPhi\si=\si,$$
where ``$\cong$'' is replaced by ``$=$'' since $\PD$ is separated. As an application of Theorems \ref{QCat-representation} and \ref{QCat-representation-dense}, we may present a much easier proof for the representation theorem of $\MPhi$ given in \cite{Shen2013a,Shen2016}:

\begin{thm} {\rm\cite{Shen2013a,Shen2016}} \label{Mphi-represnetation}
For any $\CQ$-distributor $\Phi:\CD\oto\CE$, a separated total $\CQ$-category $\CX$ is isomorphic to $\MPhi$ if, and only if, there exist a dense $\CQ$-functor $F:\CD\to\CX$ and a codense $\CQ$-functor $G:\CE\to\CX$ with
$$\Phi=G^{\nat}\circ F_{\nat},\quad\text{i.e.,}\quad\Phi=\CX(F-,G-).$$
\end{thm}

\begin{proof}
First of all, we show that
\begin{equation} \label{Phi-Y-uPhi-Y}
\Phi=(\sYd_{\CE})^{\nat}\circ(\uPhi)_{\nat}\circ(\sY_{\CD})_{\nat}.
\end{equation}
Indeed, since
\begin{align*}
\Phi(X,Y)&=\CE(Y,-)\rda(\Phi\lda\CD(-,X))\\
&=\sYd_{\CE}Y\rda\uPhi\sY_{\CD}X\\
&=\PdE(\uPhi\sY_{\CD}X,\sYd_{\CE}Y)\\
&=\PdE(-,\sYd_{\CE}Y)\circ\PdE(\uPhi\sY_{\CD}X,-)\\
&=(\sYd_{\CE})^{\nat}(-,Y)\circ(\uPhi\sY_{\CD})_{\nat}(X,-)
\end{align*}
for any $X\in\CD_0$, $Y\in\CE_0$, it follows that $\Phi=(\sYd_{\CE})^{\nat}\circ(\uPhi\sY_{\CD})_{\nat}=(\sYd_{\CE})^{\nat}\circ(\uPhi)_{\nat}\circ(\sY_{\CD})_{\nat}$.

Now, for the ``only if'' part, note that Theorem \ref{QCat-representation} ensures the existence of surjective $\CQ$-functors $F:\PD\to\CX$ and $G:\PdE\to\CX$ with
$$(\uPhi)_{\nat}=G^{\nat}\circ F_{\nat}.$$
Since $\sY_{\CD}:\CD\to\PD$ is dense and $\sYd_{\CE}:\CE\to\PdE$ is codense, it is easy to deduce the density of $F\sY_{\CD}:\CD\to\CX$ and the codensity of $G\sYd_{\CE}:\CE\to\CX$ by the surjectivity of $F$ and $G$. Moreover, by Equation \eqref{Phi-Y-uPhi-Y} we immediately have
$$\Phi=(\sYd_{\CE})^{\nat}\circ(\uPhi)_{\nat}\circ(\sY_{\CD})_{\nat}=(\sYd_{\CE})^{\nat}\circ G^{\nat}\circ F_{\nat}\circ(\sY_{\CD})_{\nat}=(G\sYd_{\CE})^{\nat}\circ(F\sY_{\CD})_{\nat}.$$
$$\bfig
\morphism(-800,0)<500,300>[\CD`\PD;\sY_{\CD}]
\morphism(800,0)<-500,300>[\CE`\PdE;\sYd_{\CE}]
\morphism(-800,0)|b|<800,-300>[\CD`\CX;F\sY_{\CD}]
\morphism(800,0)|b|<-800,-300>[\CE`\CX;G\sYd_{\CE}]
\morphism(-300,300)|a|/@{>}@<4pt>/<600,0>[\PD`\PdE;\uPhi]
\morphism(300,300)|b|/@{>}@<4pt>/<-600,0>[\PdE`\PD;\dPhi]
\morphism(-300,300)|l|<300,-600>[\PD`\CX;F]
\morphism(300,300)|r|<-300,-600>[\PdE`\CX;G]
\place(0,305)[\mbox{\scriptsize$\bot$}]
\efig$$

Conversely, for the ``if'' part, since we already have dense $\CQ$-functors $\sY_{\CD}:\CD\to\PD$, $F:\CD\to\CX$ and codense $\CQ$-functors $\sYd_{\CE}:\CE\to\PdE$, $G:\CE\to\CX$ with
$$G^{\nat}\circ F_{\nat}=\Phi=(\sYd_{\CE})^{\nat}\circ(\uPhi)_{\nat}\circ(\sY_{\CD})_{\nat}
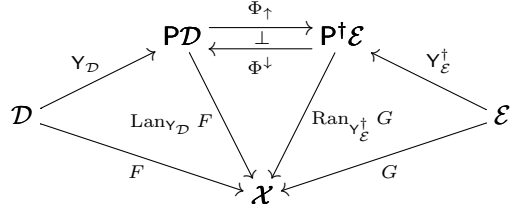
$$
$$\bfig
\morphism<600,300>[\CD`\PD;\sY_{\CD}]
\morphism(1800,0)<-600,300>[\CE`\PdE;\sYd_{\CE}]
\morphism|b|<900,-300>[\CD`\CX;F]
\morphism(1800,0)|b|<-900,-300>[\CE`\CX;G]
\morphism(600,300)|a|/@{>}@<4pt>/<600,0>[\PD`\PdE;\uPhi]
\morphism(1200,300)|b|/@{>}@<4pt>/<-600,0>[\PdE`\PD;\dPhi]
\morphism(600,300)|l|<300,-600>[\PD`\CX;\Lan_{\sY_{\CD}}F]
\morphism(1200,300)|r|<-300,-600>[\PdE`\CX;\Ran_{\sYd_{\CE}}G]
\place(900,305)[\mbox{\scriptsize$\bot$}]
\efig$$
by Equation \eqref{Phi-Y-uPhi-Y}, the conclusion follows soon from Theorem \ref{QCat-representation-dense}, which completes the proof.
\end{proof}

In the case that $\Phi$ is the identity $\CQ$-distributor $\CE:\CE\oto\CE$, $\ME$ is precisely the \emph{MacNeille completion} \cite{Shen2013a,Garner2014,Shen2014} of the $\CQ$-category $\CE$, and it can be alternatively characterized as follows by Theorem \ref{Mphi-represnetation}:

\begin{cor} \label{ME-representation}
A total $\CQ$-category $\CD$ is the MacNeille completion of its full $\CQ$-subcategory $\CE$ if $\CE$ is simultaneously dense and codense in $\CD$.
\end{cor}

The above corollary in conjunction with the definition of the MacNeille completion of a concrete category reproduces the following result of Garner:

\begin{cor} {\rm\cite{Garner2014}} \label{ME-representation-concrete}
A topological category $\CD$ over $\CB$ is the MacNeille completion of its full subcategory $\CE$ if, and only if, $\oD$ is the MacNeille completion of its $\QB$-subcategory $\oE$.
\end{cor}

\section*{Acknowledgement}

The authors acknowledge the support of National Natural Science Foundation of China (No. 11771310 and No. 11701396). We thank the anonymous referee for several helpful remarks.





\end{document}